\documentclass[final,1p,times]{elsarticle}
\usepackage{amsmath}
\usepackage{amssymb}
\usepackage{url}
\usepackage{amsfonts}
\usepackage{graphicx}
\usepackage{rotating}
\usepackage{floatflt,epsfig}
\usepackage{lineno,hyperref}
\usepackage{enumerate}
\usepackage{colortbl}
\usepackage{array,tabularx,tabulary,booktabs}
\usepackage{longtable}
\usepackage{multirow}
\usepackage{wrapfig}
\usepackage{subcaption}
\newcolumntype{^}{>{\currentrowstyle}}

\makeatletter
\def\ps@pprintTitle{%
   \let\@oddhead\@empty
   \let\@evenhead\@empty
   \let\@oddfoot\@empty
   \let\@evenfoot\@oddfoot
}
\makeatother

\newtheorem{lemma}{Lemma}

\newtheorem{construction}{Construction}
\newtheorem{property}{Property}

\begin{document}
\renewcommand{\abstractname}{Abstract}
\renewcommand{\refname}{References}
\renewcommand{\tablename}{Table}
\renewcommand{\arraystretch}{0.9}
\thispagestyle{empty}
\sloppy

\begin{frontmatter}
\title{Enumeration of strictly Deza graphs with at most 21 vertices}

\author[01]{Sergey~Goryainov}
\ead{sergey.goryainov3@gmail.com}

\author[01,02]{Dmitry~Panasenko}
\ead{makare95@mail.ru}

\author[01]{Leonid~Shalaginov}
\ead{44sh@mail.ru}

\address[01]{Chelyabinsk State University, Brat'ev Kashirinyh st. 129, Chelyabinsk, 454021, Russia}

\address[02]{Krasovskii Institute of Mathematics and Mechanics, S. Kovalevskaja st. 16, Yekaterinburg, 620990, Russia}

\begin{abstract}
A Deza graph $\Gamma$ with parameters $(v,k,b,a)$ is a $k$-regular graph with $v$ vertices such that any two distinct vertices have $b$ or $a$ common neighbours, where $b \geqslant a$. A Deza graph of diameter 2 which is not a strongly regular graph is called a strictly Deza graph. We find all 139 strictly Deza graphs up to 21 vertices and list corresponding constructions and properties. 
\end{abstract}

\begin{keyword} Deza graph \sep strictly Deza graph \sep strongly regular graph \sep dual Seidel switching
\vspace{\baselineskip}
\MSC[2010] 05C50\sep 05E10\sep 15A18
\end{keyword}
\end{frontmatter}

\section{Introduction}

Deza graphs were introduced in 1999~\cite{EFHHH1999} as a generalisation of strongly regular graphs. A \emph{Deza graph} $\Gamma$ with parameters $(v,k,b,a)$ is a $k$-regular graph with $v$ vertices for which the number of common neighbours of two distinct vertices takes just two values, $b$ or $a$, where $b \geqslant a$. A \emph{strongly regular graph} $G$ with parameters $(v,k,\lambda,\mu)$ is a $k$-regular graph with $v$ vertices such that any two adjacent vertices have $\lambda$ common neighbours and any two non-adjacent vertices have $\mu$ common neighbours. A Deza graph of diameter 2 which is not a strongly regular graph is called a \emph{strictly Deza graph}.

In 1999~\cite{EFHHH1999} the complete list of strictly Deza graphs with at most 13 vertices was presented and different constructions for those graphs were discussed. In 2011~\cite{GS2011} this list was extended up to 16 vertices. 
In 2014 S.~Goryainov and L.~Shalaginov~\cite{GS2014} found all Cayley-Deza graphs with $a > 0$ up to 59 vertices and listed all corresponding groups. These results are available on the web pages \url{http://alg.imm.uran.ru/dezagraphs/dezatab.html} and \url{http://alg.imm.uran.ru/dezagraphs/deza_cayleytab.html}.

A $k$-regular graph is called a \emph{divisible design graph} if its vertex set can be partitioned into $m$ classes of size $n$, such that two distinct vertices from the same class have exactly $\lambda_1$ common neighbors, and two vertices from different classes have exactly $\lambda_2$ common neighbors. The definition implies that all divisible design graphs are Deza graphs. Divisible design graphs were first studied in master's thesis by M.A. Meulenberg~\cite{M2008} and the list of feasible parameters of divisible design graphs up to 50 vertices was presented. In 2011-2013 in the following papers~\cite{CH2014,HKM2011} divisible design graphs were studied in more details and the existence of graphs was resolved in all but one cases for graphs up to 27 vertices. 

In this paper we find all strictly Deza graphs up to 21 vertices. It turns out that the number $Num(v)$ of non-isomorphic strictly Deza graphs with $v \leqslant 21$ is given by the following table: 

\begin{table}[h]
\label{tab:0}
\begin{tabular}{|c|rrrrrrrrrrrrrr|}
\hline $v$ & 8 & 9 & 10 & 11 & 12 & 13 & 14 & 15 & 16 & 17 & 18 & 19 & 20 & 21 \\ \hline
$Num(v)$ & 3 & 2 & 1 & 0 & 6 & 1 & 1 & 1 & 10 & 3 & 13 & 11 & 56 & 31 \\ \hline
\end{tabular}
\end{table}

In previous papers (\cite{EFHHH1999, GS2011}) corresponding constructions were given for almost all graphs. However, among graphs we found, the constructions of almost half of the graphs were unknown. Therefore, we give a list of the found graphs and indicate all known constructions. This shows for which graphs the constructions are unknown, thus giving one of the approaches to finding new constructions. 

This paper is organised as follows. In Section 2 we describe the algorithm used for enumerating Deza graphs. In Section 3 we give an overview of some known constructions and properties of Deza graphs. In Section 4 we present tables with enumeration results, corresponding properties and constructions and in Section 5 we take a closer look at Deza graphs with WL-rank 4 we found.

\section{Enumeration algorithm}
\subsection{Search for feasible parameters} 

Let $\Gamma$ be a Deza graph with parameters $(v,k,b,a)$. For a fixed vertex $u$ in $\Gamma$, define  
$$\alpha = |\{w \in V(\Gamma): |N(u) \cap N(w)| = a\}|$$ and  
$$\beta = |\{w \in V(\Gamma): |N(u) \cap N(w)| = b\}|,$$
where $V(\Gamma)$ is the vertex set of $\Gamma$ and $N(u), N(w)$ are the neighborhoods of $u$ and $w$, respectively.

In \cite[Proposition 1.1]{EFHHH1999} it was proved that $\alpha$ and $\beta$ do not depend on $u$ and can be computed as follows: 
$$\alpha = \frac{b(v - 1) - k(k - 1)}{b - a}, ~\beta = \frac{a(v - 1) - k(k - 1)}{a - b} \text{ if $a \neq b$}$$ 
and 
$$\alpha = \beta = \frac{k(k - 1)}{a} \text{ otherwise.}$$

At the first step, for a fixed number of vertices, we calculate all feasible parameters of Deza graphs satisfying restrictions given by the following lemma.

\begin{lemma}\label{l1}~{\rm \cite[Corollary 1.2]{EFHHH1999}} 
If there is a Deza graph with parameters $(v, k, b, a)$, then the following statements hold:

(i) $b - a$ divides $b(v - 1) - k(k - 1)$;

(ii) if $\alpha \neq 0$, then $v \geqslant 2k - a$;

(iii) if $\alpha, \beta \neq 0$, then $a(v - 1) < k(k - 1) < b(v - 1)$.
\end{lemma}

The number $Num'(v)$ of feasible parameters of Deza graphs with $v \leqslant 21$ is given in the table below:

\begin{table}[h]
\label{tab:00}
\begin{tabular}{|c|rrrrrrrrrrrrrr|}
\hline $v$ & 8 & 9 & 10 & 11 & 12 & 13 & 14 & 15 & 16 & 17 & 18 & 19 & 20 & 21 \\ \hline
$Num'(v)$ & 14 & 10 & 24 & 19 & 34 & 26 & 44 & 34 & 73 & 40 & 74 & 60 & 86 & 77\\ \hline
\end{tabular}
\end{table}

\subsection{Constructing adjacency matrices}

Given feasible parameters $(v, k, b, a)$ of a Deza graph $\Gamma$, we initially construct the first two rows of the adjacency matrix using the following method.

Let us consider three possible cases.

\textbf{Case 1:} $\alpha < k$ holds. Then there exist two adjacent vertices in $\Gamma$ that have exactly $b$ common neighbors. Without loss of generality, we may assume that the first two rows of the adjacency matrix of $\Gamma$ look like this:

~

$01\; \overbrace{1 \ldots 1}^{b} \; \overbrace{1 \ldots 1}^{k-b-1} \; \overbrace{0 \ldots 0}^{k-b-1} \; \overbrace{0 \ldots 0}^{v-2k+b}$

$10\;\, 1 \ldots 1 \;\, 0 \ldots 0 \;\, 1 \ldots 1 \;\, 0 \ldots 0$

~

\textbf{Case 2:} $\alpha = k$ holds. Then there exist two adjacent vertices in $\Gamma$ that have exactly $b$ common neighbors (otherwise, $\Gamma$ is strongly regular). Without loss of generality, we may assume that the first two rows of the adjacency matrix of $\Gamma$ are the same as in Case 1.

\textbf{Case 3:} $\alpha > k$ holds.
Then there exist two non-adjacent vertices in $\Gamma$ that have exactly $a$ common neighbors. Without loss of generality, we may assume that the first two rows of the adjacency matrix of $\Gamma$ look like this:

~

$00\; \overbrace{1 \ldots 1}^{a} \; \overbrace{1 \ldots 1}^{k-a} \; \overbrace{0 \ldots 0}^{k-a} \; \overbrace{0 \ldots 0}^{v-2k+a-2}$

$00\;\, 1 \ldots 1 \;\, 0 \ldots 0 \;\, 1 \ldots 1 \;\;\; 0 \ldots 0$

~

To enumerate the next three rows of the matrix, we use the following approach.

Assume we have specified the first two rows of the adjacency matrix of a Deza graph $\Gamma$:

~

\begin{tabular}{rl}
$0\spadesuit1$ & $1\ldots1~1\ldots1~0\ldots0~0\ldots0$\\
$\spadesuit01$ & $1\ldots1~0\ldots0~1\ldots1~0\ldots0$,
\end{tabular}

~\\
where $\spadesuit$ denotes the same symbol (`1' in Cases 1 and 2; `0' in Case 3).

Further, we construct the third row. Notice that the columns can be divided into four blocks as follows:

~

\begin{tabular}{rc|c|c|c}
$0\spadesuit1$ & $1\ldots1$ & $1\ldots1$ & $0\ldots0$ & $0\ldots0$\\
$\spadesuit01$ & $1\ldots1$ & $0\ldots0$ & $1\ldots1$ & $0\ldots0$\\
$110$ & $\ast\ldots\ast$ & $\ast\ldots\ast$ & $\ast\ldots\ast$ & $\ast\ldots\ast$
\end{tabular}.

~\\
and changing the order of 1s inside each block gives an equivalent matrix. We say that two partially filled matrices are \emph{equivalent} if the graphs determined by them are isomorphic. Note that, for a partially filled matrix, changing the order of vertices inside each block gives an equivalent matrix.

For example, these two matrices are equivalent:

~

\begin{tabular}{rl|l|l|lrrrrl|l|l|l}
011 & 111 & 111 & 000 & 00 & ~ & ~ & ~ & 011 & 111 & 111 & 000 & 00\\
101 & 111 & 000 & 111 & 00 & ~ & ~ & ~ & 101 & 111 & 000 & 111 & 00\\
110 & 110 & 100 & 110 & 10 & ~ & ~ & ~ & 110 & 101 & 001 & 011 & 10
\end{tabular}

~

For the next rows, this division into blocks can be extended, where the number of blocks multiplies by 2 with each row (8 blocks for the 4th row, 16 blocks for the 5th row, etc.).

Thus, to construct the next row of the adjacency matrix, we consider all possible numbers of 1s in each block. Then the obtained matrices are checked for equivalence using Magma, and the procedure for adding a new row repeats for all nonequivalent matrices. Since the equivalent matrices will give isomorphic graphs at the end, leaving all nonequivalent options does not reduce the exhaustiveness of the algorithm.

For the remaining $v-5$ rows of the matrix we use exhaustive search of possible rows: all possible combinations of required number of 1s in $v - i$ positions, where $i$ is the number of current row. For each added row we also check if the resulting matrix is the adjacency matrix of the Deza graph.

We use Magma to check whether graphs are isomorphic after the completion of the enumeration. In case of $a = 0$ we calculate the diameter of the resulting graphs. If the diameter does not equal to 2, then this graph is not a strictly Deza graph.

\section{Constructions and properties of Deza graphs}

\subsection{Cayley graphs}

Let $G$ be a group and $S \subset G$ be an identity-free subset with the property $S=S^{{-1}}$ (that is for each $s\in S \; s^{-1}$ also is an element of $S$). The \emph{Cayley graph} $Cay(G, S)$ of the group $G$ with the generating set $S$ is the graph whose vertices are elements of the group $G$ and the set of edges is given by $\{\{g,gs\}~:~g \in G, s \in S\}$.

Let $SS^{-1}$ denote the multiset $\{ss'^{-1}~:~s,s' \in S\}$ and the writing $SS^{-1} = aA + bB + k\{e\}$ mean that $SS^{-1}$ contains $a$ copies of each element of $A$, $b$ copies of each element of $B$ and $k$ copies of $e$.

\begin{construction}[{\cite[Proposition 2.1]{EFHHH1999}}]\label{c1}
A Cayley graph of a group $G$ with the generating set $S$ $Cay(G, S)$ is a Deza graph with parameters $(v, k, b, a)$ if and only if the following two conditions are satisfied:

(i) $|G| = v$ and $|S| = k$;

(ii) $SS^{-1} = aA + bB + k\{e\}$, where $A, B, \{e\}$ are a partition of $G$.
\end{construction}

In the resulting table below, we denote Cayley-Deza graphs as `cay'.

\subsection{Association schemes}

Let $X$ be a set of size $n$, and $R_0, R_1, \ldots, R_d$ be relations defined on $X$. Let $A_0, A_1, \ldots, A_d$ be the 0-1 matrices corresponding to the relations, that is, the $(x, y)$-entry of $A_i$ is 1 if and only if $(x, y) \in R_i$. Then $(X, \{R_i\}_{i=0}^d )$ is called a d-class \emph{symmetric association scheme} if

(i) $A_0 = I$, where $I$ is the identity matrix;

(ii) $\sum_{A_i} = J$, where $J$ is the all-ones matrix;

(iii) each $A_i$ is symmetric;

(iv) for each pair $i$ and $j$, $A_iA_j = \sum_k p^k_{ij} A_k$ for some constants $p^k_{ij}$.

\begin{construction}[{\cite[Theorem 4.2]{EFHHH1999}}]\label{c2}
Let $(X, \{R_0, R_1, \ldots, R_d\})$ be a symmetric association scheme, and $F \subset \{1, 2, \ldots, d\}$. Let $\Gamma$ be the graph with adjacency matrix $\sum_{f \in F} A_f$. Then $\Gamma$ is a Deza graph if and only if
$$\sum\limits_{f,g \in F} p^k_{fg}$$
takes on at most two values, as $k$ ranges over $\{1, 2, \ldots, d\}$.
\end{construction}

In the resulting table below, we denote Deza graphs obtained from association schemes as `as'.

\subsection{Dual Seidel switching}

An involutive automorphism of a graph is called \emph{Seidel automorphism} if it interchanges only non-adjacent vertices. 
\begin{construction}[{Dual Seidel switching; \cite[Theorem 3.1]{EFHHH1999}}]\label{c3}  
Let $G$ be a strongly regular graph with parameters $(v,k,\lambda,\mu)$, where $k\neq \mu$, $\lambda \neq \mu$. Let $M$ be the adjacency matrix of $G$, and $P$ be a non-identity permutation matrix of the same size. Then $PM$ is the adjacency matrix of a Deza graph $\Gamma$ if and only if $P$ represents a Seidel automorphism. Moreover, $\Gamma$ is a strictly Deza graph if and only if $\lambda \neq 0$, $\mu \neq 0$.
\end{construction}

\begin{construction}[{Generalised dual Seidel switching 1; \cite[Theorem 5]{KKS2021}}]\label{c4}
Let $G$ be a strongly regular graph with the adjacency matrix $M$, and $H$ be its induced subgraph with the adjacency matrix $M_{11}$. If there exists a Seidel automorphism of $H$ with the permutation matrix $P_{11}$ such that $P_{11}M_{12}M_{22}=M_{12}M_{22}$, then matrix
$$N= \left(\begin{array}{cc} P_{11}M_{11} & M_{12}\\ M_{21}& M_{22}\\ \end{array}\right)$$
is the adjacency matrix of a Deza graph.
\end{construction}

In the resulting table below, we denote Deza graphs obtained by dual Seidel switching as `dss' and Deza graphs obtained by generalised dual Seidel switching as `gdss'.

\begin{construction}[{Generalised dual Seidel switching 2; \cite[Theorem 6]{KKS2021}}]\label{c4.1}
Let $\Gamma$ be a Deza graph with the adjacency matrix $M$, and $H$ be its induced subgraph with the adjacency matrix $M_{11}$. If there exists a Seidel automorphism of $H$ with the permutation matrix $P_{11}$ such that $P_{11}M_{12}M_{22}=M_{12}M_{22}$, then matrix
$$N= \left(\begin{array}{cc} P_{11}M_{11} & M_{12}\\ M_{21}& M_{22}\\ \end{array}\right)$$
is the adjacency matrix of a Deza graph.
\end{construction}

Note that in \cite{KKS2021} this construction was considered only for Deza graphs with strongly regular children (see definitions in Section 3.5) but the proof does not use this property, therefore this construction can be applied to any Deza graph.

In the resulting table below, we denote Deza graphs obtained by generalised dual Seidel switching from Deza graphs as `gdss($n$)', where $n$ denotes the serial number of the used Deza graph from the table.

\begin{construction}[{\cite[Theorem 7]{KKS2021}}]\label{c5}
Let $M$ be the adjacency matrix of a strongly regular graph $G$ with parameters $(v,k,\lambda,\mu)$ with $\lambda=\mu$. If there exists a fixed point free Seidel automorphism of $G$ and its permutation matrix is $P$, then the matrix $M+P$ is the adjacency matrix of a Deza graph.
\end{construction}

\begin{construction}[{\cite[Theorem 8]{KKS2021}}]\label{c6}
Let $M$ be the adjacency matrix of a strongly regular graph $G$ with parameters $(v,k,\lambda,\mu)$. If there exists a fixed point free Seidel automorphism of $G$ and its permutation matrix is $P$, then the matrix $P(M+I)$ is the adjacency matrix of a Deza graph.
\end{construction}

In the resulting table below, we denote Deza graphs obtained by construction 6 as `c6' and Deza graphs obtained by construction 7 as `c7'.

\subsection{Lexicographic product of graphs}
The \emph{lexicographic product} or \emph{graph composition} $G[H]$ of graphs $G$ and $H$ is a graph such that the vertex set of $G[H]$ is $V(G) \times V(H)$ and adjacency defined by 
$$(u_1, u_2) \sim (v_1, v_2) \text{ if and only if } u_1 \sim v_1 \text{ or } (u_1 = v_1 \text{ and } u_2 \sim v_2).$$

\begin{construction}[{\cite[Proposition 2.3]{EFHHH1999}}]\label{c7}
Let $G$ be a strongly regular graph with parameters $(v, k, \lambda, \mu)$ and $\Gamma$ be a Deza graph with parameters $(v', k', b, a)$. Then $G[\Gamma]$ is a $(k' + kv')$-regular graph on $vv'$ vertices. It is a Deza graph if and only if 
$$|\{a + kv', b + kv', \mu v', \lambda v' + 2k'\}| \leqslant 2.$$
\end{construction}

In this paper we restrict ourselves to the following three applications of this construction.

~
 
\noindent {\bf Construction~8.1.}
\emph{Let $G$ be $K_{x \times y}$, the complete multipartite graph containing $x$ parts of $y$ vertices. Then $G[K_2]$ is a Deza graph with parameters $(2xy, 2y(x - 1) + 1, 2y(x - 1), 2y(x - 2) + 2)$.} 

~

\noindent {\bf Construction~8.2.}
\emph{Let $G$ be a strongly regular graph with parameters $(v, k, \lambda, \mu)$, where $\lambda = \mu - 1$. Then $G[K_2]$ is a Deza graph with parameters $(2v, 2k + 1, 2k, 2\mu)$.} 

~

\noindent {\bf Construction~8.3.}
\emph{Let $\Gamma$ be a Deza graph obtained with construction 8.2. Let $M$ be the adjacency matrix of $\Gamma$, and $P$ be a non-identity permutation matrix of the same size. Then $PM$ is the adjacency matrix of a Deza graph if and only if $P$ represents a Seidel automorphism.} 

~

These three constructions were considered in two papers \cite{GHKS2019, KMS2019} on Deza graphs with parameters $(v,k,k-1,a)$.

In the resulting table below, we denote Deza graphs obtained by these constructions as `c8.1', `c8.2' and `c8.3'.

\subsection{Properties of Deza graphs}
Suppose $\Gamma$ is a graph with $v$ vertices, and $M$ is its adjacency matrix. Then $\Gamma$ is a Deza graph with parameters $(v,k,b,a)$ if and only if 
$$M^2 = aA + bB + kI$$
for some symmetric $(0,1)$-matrices $A$, $B$ such that $A+B+I=J$ \cite{EFHHH1999}. Note  that $\Gamma$ is a strongly regular graph if and only if $A$ or $B$ is $M$. 

Suppose that we have a Deza graph with $M$, $A$, and $B$ satisfying the equality above. Then $A$ and $B$ are the adjacency matrices of graphs, and the corresponding graphs $\Gamma_A$ and $\Gamma_B$ are called the \emph{children} of $\Gamma$.

The definition of divisible design graphs implies the following property, which can be used to determine if a Deza graph is a divisible design graph.

\begin{property}\label{p1}
A Deza graph whose children are a complete multipartite graph and a union of complete graphs is a divisible design graph. 
\end{property}

A \emph{coherent configuration} $\mathcal{X}$ on a finite set $V$ can be thought as a special partition of $V \times V$ for which the diagonal of $V \times V$ is a union of classes \cite{CP2019}. If in a coherent configuration the diagonal of $V \times V$ is a single class then this coherent configuration is an association scheme.

Each graph has a specific coherent configuration associated with it, known as \emph{WL-closure}, which can be obtained using \emph{Weisfeiler-Leman algorithm} \cite{LW1968}. Given a graph $G$ with the vertex set $V(G)$ and the edge set $E(G)$, this algorithm constructs the smallest coherent configuration on $V(G)$ for which $E(G)$ is a union of classes. The number of classes in WL-closure is called \emph{WL-rank}. 

A graph has WL-rank~3 if and only if this graph is a strongly regular graph \cite[Section 2.6.3]{CP2019}. So it is interesting to study graphs with small WL-rank more than~3. In Section~5, we take a closer look at the Deza graphs with WL-rank~4 we found. This section, in particular, shows that all such graphs we found can be described using known constructions. 

\section{Enumeration results}

In the table below, \# gives a serial number, $v, k, b, a$ are the parameters of a Deza graph, `egv' denotes the number of distinct eigenvalues, `int' denotes whether a graph has an integral spectrum, `ddg' means divisible design graph, `WL-rank' denotes WL-rank of a graph. The column `constructions' describes constructions from Section~3, which can be used to obtain this graph.

Note that sometimes generalised dual Seidel switching produces Deza graphs from Deza graphs with an unknown construction. For example, we can obtain the graph with serial number 29 from the graph with serial number 30 and vice versa. These cases are not presented in the resulting table.

~

\begin{longtable}{|c|c|c|c|c|c|c|c|c|l|}

\caption{Strictly Deza graphs with at most 21 vertices}
\label{tab:1}\\

\multicolumn{1}{|c|}{\#} & 
\multicolumn{1}{c|}{$v$} & 
\multicolumn{1}{c|}{$k$} & 
\multicolumn{1}{c|}{$b$} & 
\multicolumn{1}{c|}{$a$} & 
\multicolumn{1}{c|}{egv} & 
\multicolumn{1}{c|}{int} & 
\multicolumn{1}{c|}{ddg} & 
\multicolumn{1}{c|}{WL-rank} & 
\multicolumn{1}{c|}{constructions} \\ [1pt] \hline 
& & & & & & & & & \\ [-0.7em]
\endfirsthead

\multicolumn{1}{|c|}{\#} & 
\multicolumn{1}{c|}{$v$} & 
\multicolumn{1}{c|}{$k$} & 
\multicolumn{1}{c|}{$b$} & 
\multicolumn{1}{c|}{$a$} & 
\multicolumn{1}{c|}{egv} & 
\multicolumn{1}{c|}{int} & 
\multicolumn{1}{c|}{ddg} & 
\multicolumn{1}{c|}{WL-rank} & 
\multicolumn{1}{c|}{constructions} \\ [1pt] \hline 
& & & & & & & & & \\ [-0.7em]
\endhead

\hline \multicolumn{10}{c}{} \\
\endfoot

\hline \multicolumn{10}{c}{} \\
\endlastfoot

1 & 8 & 4 & 2 & 0 & 4 & + & + & 4 & cay, as \\ [2pt]
2 & 8 & 4 & 2 & 1 & 5 & -- & -- & 5 & cay, as \\ [2pt]
3 & 8 & 5 & 4 & 2 & 4 & + & + & 4 & cay, as, c8.1 \\ [2pt]
4 & 9 & 4 & 2 & 1 & 5 & + & -- & 10 & dss, gdss \\ [2pt]
5 & 9 & 4 & 2 & 1 & 5 & -- & -- & 5 & cay, as \\ [2pt]
6 & 10 & 5 & 4 & 2 & 4 & -- & + & 4 & cay, as, c8.2 \\ [2pt]
7 & 12 & 5 & 2 & 1 & 4 & + & + & 4 & cay, as \\ [2pt]
8 & 12 & 6 & 3 & 2 & 5 & -- & -- & 10 & gdss(9) \\ [2pt]
9 & 12 & 6 & 3 & 2 & 4 & + & + & 5 & cay, as \\ [2pt]
10 & 12 & 7 & 4 & 3 & 5 & + & + & 6 & cay, as \\ [2pt]
11 & 12 & 7 & 6 & 2 & 4 & + & + & 4 & cay, as, c8.1 \\ [2pt]
12 & 12 & 9 & 8 & 6 & 4 & + & + & 4 & cay, as, c8.1 \\ [2pt]
13 & 13 & 8 & 5 & 4 & 4 & -- & -- & 4 & cay, as \\ [2pt]
14 & 14 & 9 & 6 & 4 & 4 & -- & -- & 4 & cay, as \\ [2pt]
15 & 15 & 6 & 3 & 1 & 5 & + & -- & 16 & dss, gdss \\ [2pt]
16 & 16 & 5 & 2 & 1 & 7 & -- & -- & 16 & cay \\ [2pt]
17 & 16 & 7 & 4 & 2 & 5 & + & -- & 6 & cay, as, c6, c7 \\ [2pt]
18 & 16 & 7 & 4 & 2 & 5 & + & -- & 6 & cay, as, c6, c7 \\ [2pt]
19 & 16 & 8 & 4 & 2 & 6 & -- & -- & 8 & cay, as \\ [2pt]
20 & 16 & 9 & 6 & 4 & 5 & + & -- & 6 & cay, as, dss, gdss \\ [2pt]
21 & 16 & 9 & 6 & 4 & 5 & + & -- & 12 & dss, gdss \\ [2pt]
22 & 16 & 9 & 8 & 2 & 4 & + & + & 4 & cay, as, c8.1 \\ [2pt]
23 & 16 & 11 & 8 & 6 & 5 & + & -- & 5 & cay, as, c6, c7 \\ [2pt]
24 & 16 & 12 & 10 & 8 & 5 & + & -- & 5 & cay, as \\ [2pt]
25 & 16 & 13 & 12 & 10 & 4 & + & + & 4 & cay, as, c8.1 \\ [2pt]
26 & 17 & 8 & 4 & 3 & 10 & -- & -- & 93 & -- \\ [2pt]
27 & 17 & 8 & 4 & 3 & 13 & -- & -- & 83 & -- \\ [2pt]
28 & 17 & 8 & 4 & 3 & 13 & -- & -- & 83 & -- \\ [2pt]
29 & 18 & 8 & 4 & 2 & 18 & -- & -- & 162 & -- \\ [2pt]
30 & 18 & 8 & 4 & 2 & 12 & -- & -- & 34 & -- \\ [2pt]
31 & 18 & 8 & 4 & 2 & 13 & -- & -- & 65 & -- \\ [2pt]
32 & 18 & 8 & 4 & 2 & 10 & -- & -- & 18 & cay \\ [2pt]
33 & 18 & 8 & 4 & 2 & 11 & -- & -- & 54 & -- \\ [2pt]
34 & 18 & 8 & 4 & 2 & 8 & -- & -- & 19 & gdss(32, 35) \\ [2pt]
35 & 18 & 8 & 4 & 2 & 5 & -- & -- & 5 & cay, as \\ [2pt]
36 & 18 & 8 & 4 & 3 & 13 & -- & -- & 98 & -- \\ [2pt]
37 & 18 & 9 & 6 & 4 & 7 & -- & -- & 36 & gdss(38) \\ [2pt]
38 & 18 & 9 & 6 & 4 & 5 & -- & + & 5 & cay, as \\ [2pt]
39 & 18 & 9 & 8 & 4 & 5 & + & + & 13 & c8.3 \\ [2pt]
40 & 18 & 9 & 8 & 4 & 4 & + & + & 4 & cay, as, c8.2 \\ [2pt]
41 & 18 & 13 & 12 & 8 & 4 & + & + & 4 & cay, as, c8.1 \\ [2pt]
42 & 19 & 6 & 2 & 1 & 13 & -- & -- & 65 & -- \\ [2pt]
43 & 19 & 6 & 2 & 1 & 13 & -- & -- & 65 & -- \\ [2pt]
44 & 19 & 6 & 2 & 1 & 13 & -- & -- & 65 & -- \\ [2pt]
45 & 19 & 6 & 2 & 1 & 13 & -- & -- & 65 & -- \\ [2pt]
46 & 19 & 6 & 2 & 1 & 4 & -- & -- & 4 & cay, as \\ [2pt]
47 & 19 & 6 & 2 & 1 & 13 & -- & -- & 65 & -- \\ [2pt]
48 & 19 & 8 & 4 & 2 & 9 & -- & -- & 55 & -- \\ [2pt]
49 & 19 & 8 & 4 & 2 & 14 & -- & -- & 93 & -- \\ [2pt]
50 & 19 & 8 & 4 & 2 & 18 & -- & -- & 361 & -- \\ [2pt]
51 & 19 & 12 & 8 & 7 & 8 & -- & -- & 24 & -- \\ [2pt]
52 & 19 & 12 & 8 & 7 & 13 & -- & -- & 61 & -- \\ [2pt]
53 & 20 & 6 & 2 & 1 & 11 & -- & -- & 42 & -- \\ [2pt]
54 & 20 & 6 & 2 & 1 & 10 & -- & -- & 100 & -- \\ [2pt]
55 & 20 & 6 & 2 & 1 & 18 & -- & -- & 200 & -- \\ [2pt]
56 & 20 & 6 & 2 & 1 & 5 & -- & -- & 6 & as \\ [2pt]
57 & 20 & 6 & 2 & 1 & 10 & -- & -- & 80 & -- \\ [2pt]
58 & 20 & 6 & 2 & 1 & 20 & -- & -- & 400 & -- \\ [2pt]
59 & 20 & 6 & 2 & 1 & 20 & -- & -- & 400 & -- \\ [2pt]
60 & 20 & 6 & 2 & 1 & 19 & -- & -- & 400 & -- \\ [2pt]
61 & 20 & 6 & 2 & 1 & 19 & -- & -- & 200 & -- \\ [2pt]
62 & 20 & 6 & 2 & 1 & 18 & -- & -- & 200 & -- \\ [2pt]
63 & 20 & 6 & 2 & 1 & 20 & -- & -- & 202 & -- \\ [2pt]
64 & 20 & 6 & 2 & 1 & 20 & -- & -- & 400 & -- \\ [2pt]
65 & 20 & 6 & 2 & 1 & 19 & -- & -- & 400 & -- \\ [2pt]
66 & 20 & 6 & 2 & 1 & 16 & -- & -- & 202 & -- \\ [2pt]
67 & 20 & 6 & 2 & 1 & 20 & -- & -- & 400 & -- \\ [2pt]
68 & 20 & 6 & 2 & 1 & 20 & -- & -- & 202 & -- \\ [2pt]
69 & 20 & 6 & 2 & 1 & 20 & -- & -- & 400 & -- \\ [2pt]
70 & 20 & 6 & 2 & 1 & 18 & -- & -- & 200 & -- \\ [2pt]
71 & 20 & 6 & 2 & 1 & 18 & -- & -- & 400 & -- \\ [2pt]
72 & 20 & 6 & 2 & 1 & 20 & -- & -- & 400 & -- \\ [2pt]
73 & 20 & 6 & 2 & 1 & 20 & -- & -- & 202 & -- \\ [2pt]
74 & 20 & 6 & 2 & 1 & 8 & -- & -- & 42 & -- \\ [2pt]
75 & 20 & 6 & 2 & 1 & 16 & -- & -- & 202 & -- \\ [2pt]
76 & 20 & 6 & 2 & 1 & 20 & -- & -- & 200 & -- \\ [2pt]
77 & 20 & 6 & 2 & 1 & 20 & -- & -- & 200 & -- \\ [2pt]
78 & 20 & 6 & 2 & 1 & 8 & -- & -- & 42 & gdss(56) \\ [2pt]
79 & 20 & 6 & 2 & 1 & 7 & -- & -- & 31 & gdss(80) \\ [2pt]
80 & 20 & 6 & 2 & 1 & 7 & -- & -- & 27 & gdss(56) \\ [2pt]
81 & 20 & 6 & 2 & 1 & 20 & -- & -- & 202 & -- \\ [2pt]
82 & 20 & 6 & 2 & 1 & 20 & -- & -- & 400 & -- \\ [2pt]
83 & 20 & 6 & 2 & 1 & 18 & -- & -- & 200 & -- \\ [2pt]
84 & 20 & 6 & 2 & 1 & 18 & -- & -- & 200 & -- \\ [2pt]
85 & 20 & 6 & 2 & 1 & 20 & -- & -- & 200 & -- \\ [2pt]
86 & 20 & 6 & 2 & 1 & 10 & -- & -- & 80 & -- \\ [2pt]
87 & 20 & 6 & 2 & 1 & 20 & -- & -- & 202 & -- \\ [2pt]
88 & 20 & 6 & 2 & 1 & 10 & -- & -- & 40 & -- \\ [2pt]
89 & 20 & 6 & 2 & 1 & 11 & -- & -- & 40 & -- \\ [2pt]
90 & 20 & 7 & 3 & 2 & 4 & + & + & 4 & cay, as \\ [2pt]
91 & 20 & 7 & 6 & 2 & 4 & + & + & 4 & cay, as, c8.2 \\ [2pt]
92 & 20 & 7 & 6 & 2 & 5 & + & + & 19 & c8.3 \\ [2pt]
93 & 20 & 8 & 4 & 2 & 15 & -- & -- & 122 & -- \\ [2pt]
94 & 20 & 8 & 4 & 2 & 14 & -- & -- & 59 & -- \\ [2pt]
95 & 20 & 8 & 4 & 2 & 14 & -- & -- & 100 & -- \\ [2pt]
96 & 20 & 8 & 4 & 2 & 5 & -- & -- & 20 & cay \\ [2pt]
97 & 20 & 8 & 4 & 2 & 13 & -- & -- & 208 & -- \\ [2pt]
98 & 20 & 10 & 6 & 4 & 5 & -- & -- & 40 & -- \\ [2pt]
99 & 20 & 10 & 6 & 4 & 13 & -- & -- & 200 & -- \\ [2pt]
100 & 20 & 10 & 6 & 4 & 13 & -- & -- & 200 & -- \\ [2pt]
101 & 20 & 10 & 6 & 4 & 9 & -- & -- & 36 & -- \\ [2pt]
102 & 20 & 10 & 6 & 4 & 5 & -- & -- & 7 & cay, as \\ [2pt]
103 & 20 & 10 & 6 & 4 & 7 & -- & -- & 40 & -- \\ [2pt]
104 & 20 & 11 & 10 & 2 & 4 & + & + & 4 & cay, as, c8.1 \\ [2pt]
105 & 20 & 13 & 9 & 8 & 5 & + & + & 6 & cay, as \\ [2pt]
106 & 20 & 13 & 12 & 8 & 4 & + & + & 4 & cay, as, c8.2 \\ [2pt]
107 & 20 & 14 & 10 & 9 & 5 & -- & -- & 6 & as \\ [2pt]
108 & 20 & 17 & 16 & 14 & 4 & + & + & 4 & cay, as, c8.1 \\ [2pt]
109 & 21 & 8 & 3 & 2 & 8 & -- & -- & 50 & -- \\ [2pt]
110 & 21 & 8 & 3 & 2 & 8 & -- & -- & 35 & -- \\ [2pt]
111 & 21 & 8 & 3 & 2 & 5 & -- & -- & 12 & cay, as \\ [2pt]
112 & 21 & 8 & 4 & 2 & 7 & -- & -- & 129 & gdss(122) \\ [2pt]
113 & 21 & 8 & 4 & 2 & 20 & -- & -- & 225 & -- \\ [2pt]
114 & 21 & 8 & 4 & 2 & 7 & -- & -- & 28 & gdss(130) \\ [2pt]
115 & 21 & 8 & 4 & 2 & 21 & -- & -- & 225 & -- \\ [2pt]
116 & 21 & 8 & 4 & 2 & 15 & -- & -- & 441 & gdss(122) \\ [2pt]
117 & 21 & 8 & 4 & 2 & 11 & -- & -- & 117 & gdss(122) \\ [2pt]
118 & 21 & 8 & 4 & 2 & 21 & -- & -- & 225 & -- \\ [2pt]
119 & 21 & 8 & 4 & 2 & 19 & -- & -- & 225 & gdss(123) \\ [2pt]
120 & 21 & 8 & 4 & 2 & 21 & -- & -- & 225 & -- \\ [2pt]
121 & 21 & 8 & 4 & 2 & 7 & -- & -- & 38 & gdss(130) \\ [2pt]
122 & 21 & 8 & 4 & 2 & 7 & -- & -- & 71 & gdss(130) \\ [2pt]
123 & 21 & 8 & 4 & 2 & 8 & -- & -- & 63 & gdss(122) \\ [2pt]
124 & 21 & 8 & 4 & 2 & 19 & -- & -- & 225 & -- \\ [2pt]
125 & 21 & 8 & 4 & 2 & 12 & -- & -- & 78 & -- \\ [2pt]
126 & 21 & 8 & 4 & 2 & 12 & -- & -- & 27 & -- \\ [2pt]
127 & 21 & 8 & 4 & 2 & 12 & -- & -- & 30 & -- \\ [2pt]
128 & 21 & 8 & 4 & 2 & 6 & -- & -- & 26 & -- \\ [2pt]
129 & 21 & 8 & 4 & 2 & 12 & -- & -- & 27 & -- \\ [2pt]
130 & 21 & 8 & 4 & 2 & 4 & -- & -- & 4 & cay, as \\ [2pt]
131 & 21 & 10 & 5 & 4 & 8 & -- & -- & 32 & -- \\ [2pt]
132 & 21 & 10 & 5 & 4 & 11 & -- & -- & 63 & -- \\ [2pt]
133 & 21 & 10 & 5 & 4 & 5 & + & -- & 46 & gdss \\ [2pt]
134 & 21 & 10 & 5 & 4 & 5 & + & -- & 117 & gdss(133) \\ [2pt]
135 & 21 & 10 & 6 & 3 & 5 & + & -- & 16 & dss, gdss \\ [2pt]
136 & 21 & 12 & 7 & 5 & 7 & -- & -- & 46 & gdss(137) \\ [2pt]
137 & 21 & 12 & 7 & 5 & 4 & -- & -- & 4 & cay, as \\ [2pt]
138 & 21 & 12 & 7 & 6 & 5 & -- & -- & 5 & cay, as \\ [2pt]
139 & 21 & 12 & 7 & 6 & 5 & -- & -- & 12 & cay, as \\ [2pt]

\end{longtable}

Among 139 Deza graphs we found there are 30 graphs with integral spectra. These graphs and their spectra are listed in the table bellow.

~

\begin{longtable}{|c|c|c|c|c|llll|}

\caption{Strictly Deza graphs with integral spectra}
\label{tab:2}\\ 

\multicolumn{1}{|c|}{\#} & 
\multicolumn{1}{c|}{$v$} & 
\multicolumn{1}{c|}{$k$} & 
\multicolumn{1}{c|}{$b$} & 
\multicolumn{1}{c|}{$a$} & 
\multicolumn{4}{c|}{non-principal eigenvalues} \\ [1pt] \hline 
& & & & & & & & \\ [-0.7em]
\endfirsthead

\multicolumn{1}{|c|}{\#} & 
\multicolumn{1}{c|}{$v$} & 
\multicolumn{1}{c|}{$k$} & 
\multicolumn{1}{c|}{$b$} & 
\multicolumn{1}{c|}{$a$} & 
\multicolumn{4}{c|}{non-principal eigenvalues} \\ [1pt] \hline 
& & & & & & & & \\ [-0.7em]
\endhead

\hline \multicolumn{9}{c}{} \\
\endfoot

\hline \multicolumn{9}{c}{} \\
\endlastfoot

1 & 8 & 4 & 2 & 0 & $-2^{3}$ & $0^{3}$ & $2^{1}$ & \\ [2pt]
3 & 8 & 5 & 4 & 2 & $-3^{1}$ & $-1^{4}$ & $1^{2}$ & \\ [2pt]
4 & 9 & 4 & 2 & 1 & $-2^{3}$ & $-1^{2}$ & $1^{2}$ & $2^{1}$ \\ [2pt]
7 & 12 & 5 & 2 & 1 & $-2^{6}$ & $1^{3}$ & $2^{2}$ & \\ [2pt]
9 & 12 & 6 & 3 & 2 & $-2^{6}$ & $0^{2}$ & $2^{3}$ & \\ [2pt]
10 & 12 & 7 & 4 & 3 & $-2^{6}$ & $-1^{1}$ & $1^{2}$ & $2^{2}$ \\ [2pt]
11 & 12 & 7 & 6 & 2 & $-5^{1}$ & $-1^{6}$ & $1^{4}$ & \\ [2pt]
12 & 12 & 9 & 8 & 6 & $-3^{2}$ & $-1^{6}$ & $1^{3}$ & \\ [2pt]
15 & 15 & 6 & 3 & 1 & $-3^{4}$ & $-1^{3}$ & $1^{6}$ & $3^{1}$ \\ [2pt]
17 & 16 & 7 & 4 & 2 & $-3^{4}$ & $-1^{5}$ & $1^{4}$ & $3^{2}$ \\ [2pt]
18 & 16 & 7 & 4 & 2 & $-3^{4}$ & $-1^{5}$ & $1^{4}$ & $3^{2}$ \\ [2pt]
20 & 16 & 9 & 6 & 4 & $-3^{4}$ & $-1^{6}$ & $1^{3}$ & $3^{2}$ \\ [2pt]
21 & 16 & 9 & 6 & 4 & $-3^{5}$ & $-1^{3}$ & $1^{6}$ & $3^{1}$ \\ [2pt]
22 & 16 & 9 & 8 & 2 & $-7^{1}$ & $-1^{8}$ & $1^{6}$ & \\ [2pt]
23 & 16 & 11 & 8 & 6 & $-3^{4}$ & $-1^{6}$ & $1^{4}$ & $3^{1}$ \\ [2pt]
24 & 16 & 12 & 10 & 8 & $-4^{1}$ & $-2^{6}$ & $0^{6}$ & $2^{2}$ \\ [2pt]
25 & 16 & 13 & 12 & 10 & $-3^{3}$ & $-1^{8}$ & $1^{4}$ & \\ [2pt]
39 & 18 & 9 & 8 & 4 & $-3^{5}$ & $-1^{6}$ & $1^{3}$ & $3^{3}$ \\ [2pt]
40 & 18 & 9 & 8 & 4 & $-3^{4}$ & $-1^{9}$ & $3^{4}$ & \\ [2pt]
41 & 18 & 13 & 12 & 8 & $-5^{2}$ & $-1^{9}$ & $1^{6}$ & \\ [2pt]
90 & 20 & 7 & 3 & 2 & $-2^{12}$ & $2^{4}$ & $3^{3}$ & \\ [2pt]
91 & 20 & 7 & 6 & 2 & $-3^{4}$ & $-1^{10}$ & $3^{5}$ & \\ [2pt]
92 & 20 & 7 & 6 & 2 & $-3^{5}$ & $-1^{7}$ & $1^{3}$ & $3^{4}$ \\ [2pt]
104 & 20 & 11 & 10 & 2 & $-9^{1}$ & $-1^{10}$ & $1^{8}$ & \\ [2pt]
105 & 20 & 13 & 9 & 8 & $-3^{1}$ & $-2^{12}$ & $2^{4}$ & $3^{2}$ \\ [2pt]
106 & 20 & 13 & 12 & 8 & $-3^{5}$ & $-1^{10}$ & $3^{4}$ & \\ [2pt]
108 & 20 & 17 & 16 & 14 & $-3^{4}$ & $-1^{10}$ & $1^{5}$ & \\ [2pt]
133 & 21 & 10 & 5 & 4 & $-3^{2}$ & $-2^{11}$ & $2^{3}$ & $3^{4}$ \\ [2pt]
134 & 21 & 10 & 5 & 4 & $-3^{4}$ & $-2^{8}$ & $2^{6}$ & $3^{2}$ \\ [2pt]
135 & 21 & 10 & 6 & 3 & $-4^{5}$ & $-1^{4}$ & $1^{10}$ & $4^{1}$ \\ [2pt]

\end{longtable}

\section{Deza graphs with WL-rank 4}

Since WL-closure is a coherent configuration in which the edges of a graph are a union of classes, all graphs with WL-rank 4 can be obtained from 3-class association schemes. These schemes were studied by Edwin R. van Dam in 1999 \cite{D1999}.

\subsection{Product construction from strongly regular graphs}

If $G$ is a strongly regular graph, then the graph $G \otimes J_n$, defined by its adjacency matrix $M \otimes J_n$, where M is the adjacency matrix of $G$  and $J_n$ is all-ones matrix of size $n$, generates 3-class association scheme (the other relations are $\overline{G} \otimes J_n$ and a disjoint union of $n$-cliques).

Deza graphs obtained by this method were described in construction 8.1 and 8.2.

\subsection{Rectangular schemes}

The \emph{rectangular scheme} $R(m, n)$ has as vertices the ordered pairs $(i, j)$, with $i = 1, \ldots , m$ and $j = 1, \ldots , n$. For two distinct pairs we can have the following three cases. They have the same first coordinate, or the second coordinate, or both coordinates are different, and the relations are defined accordingly. 

Let us consider a rectangular scheme $R(4, n)$. If we merge classes corresponding to the first and the second case, we obtain a Deza graph with parameters $(4n, n + 2, n - 2, 2)$. In the resulting table, this construction gives graphs with serial numbers 1 ($n = 2$), 7 ($n = 3$), 90 ($n = 5$). Note that these graphs are isomorphic to $4 \times n$-lattice, and in case of $n = 4$ the resulting graph is strongly regular therefore it was excluded from the table. 

\subsection{Cyclotomic schemes}

If $v$ is a prime power and $v \equiv 1 \text{ (mod  3)}$, we can define the 3-class \emph{cyclotomic association scheme} $Cycl(v)$ as follows. Let $\gamma$ be a primitive element of $GF(v)$. We take the elements of $GF(v)$ as vertices. Two vertices belong to $i$-th relation, where $i = 1, 2, 3$, if their difference equals $\gamma^{3t+i}$ for some $t$.

In the resulting table, this construction gives graphs with serial numbers 13~($v~=~13$), 46~($v~=~19)$.

\subsection{Distance-regular graphs}

A \emph{distance-regular graph} is a connected graph for which the distance relations (i.e., a pair of vertices is in $R_i$ if their distance in the graph is $i$) form an association scheme.

The \emph{Heawood graph} is the incidence graph of the Fano plane. The Heawood graph is a distance-regular graph with diameter 3. The \emph{line graph} of a graph is a graph whose vertices are the edges of the original graph and they are adjacent if they have a common vertex in the original graph. The line graph of the Heawood graph is also a distance-regular graph with diameter 3.

In the resulting table, graph with serial number 14 can be obtained from the Heawood graph and its edges are the union of the relations ``to be at the distance 1'' and ``to be at the distance 2''. Graphs with serial numbers 130 and 137 can be obtained from the line graph of the Heawood graph. Their edges are the relation ``to be at the distance 2'' or the union of the relations ``to be at the distance 1'' and ``to be at the distance 3'', respectively.

\section{Conclusion}

The complete list of strictly Deza graphs up to 21 vertices is available by \url{http://alg.imm.uran.ru/dezagraphs/dezatab.html}. This web page provides access to adjacency matrices, WL-closures and spectra of the graphs we found.

\section*{Acknowledgments}

The reported study was funded by RFBR according to the research project 20-51-53023.

\end{document}